\documentclass[11pt,a4paper]{article}

\usepackage[margin=1in]{geometry}
\usepackage{amsmath, amssymb}

\usepackage{hyperref}
\hypersetup{
    colorlinks=true,
    linkcolor=blue,
    filecolor=magenta,      
    urlcolor=blue,     
    citecolor=blue,    
}

\title{\textbf{From Objective Compromise to Structure-Enabled Intelligence: Reshaping the Scientific Paradigm of AI}}
\author{
    CHEN Shengchang \\
    School of Earth Sciences, Zhejiang University, Hangzhou, China \\
    \texttt{chenshengc@zju.edu.cn}
}
\date{}

\begin{document}

\maketitle

\begin{abstract}
Mainstream artificial intelligence systems, particularly Large Language Models (LLMs), are fundamentally constructed upon the Optimization-Driven (OD) paradigm. This paradigm can be characterized as an ill-posed optimization problem within a high-dimensional parameter space. When confronted with real-world scenarios characterized by noise and conflicting constraints, the logic of global error minimization inherently induces a strong smoothing effect. Consequently, this yields statistical compromise solutions that flatten local features, making it arduous to overcome the inherent bottleneck of Out-of-Distribution (OOD) generalization. Existing AI models generate outputs purely based on statistical probabilities. They lack the capacity to autonomously discern truth from falsehood or establish causal logic, thereby lacking genuine cognitive intelligence. To address the ontological limitations of the OD paradigm, this paper proposes a novel paradigm termed \textit{Structure-Enabled Intelligence}. This approach abandons the conventional path of relying solely on computing power scaling and data fitting to achieve emergent intelligence. Instead, it embeds mathematical structures and physical mechanisms as underlying priors into the algorithmic architecture, constraining the model's solving process with structured scientific rules. Through this paradigm shift from ``data-driven objective fitting'' to ``mechanism-structure dominance,'' the proposed approach fundamentally mitigates the issue of statistical compromise, enhances the model's extrapolation and generalization capabilities, and simultaneously reshapes human-AI synergy by upholding human subjectivity in scientific cognition.

\vspace{1em}
\noindent \textbf{Keywords:} Artificial Intelligence; Large Language Models; Optimization-Driven; Statistical Compromise; Structure-Enabled Intelligence; Paradigm Shift
\end{abstract}

\section{Introduction}
In recent years, ``AI for Science'' has emerged as a prominent research frontier across intelligent science and various scientific disciplines, with data-driven AI methods being extensively applied in scientific computing scenarios. However, fundamentally, the core training logic of the vast majority of current deep learning models remains anchored within the Optimization-Driven (OD) paradigm: relying on the Backpropagation (BP) algorithm \cite{rumelhart1986}, targeting the minimization of an objective function, iteratively searching within a high-dimensional parameter space, and matching task distributions through massive statistical data fitting.

While this paradigm exhibits excellent engineering efficacy in in-distribution data fitting and standardized task evaluations, it is essentially a non-linear regression over scattered data points, lacking physical mechanism constraints and internal logical support. When such purely statistical-driven models are deployed in tasks requiring scientific exploration, complex inversion, or reasoning in unknown scenarios, their black-box characteristics, statistical biases, and insufficient extrapolation capabilities are fully exposed. Increasing computing power, scaling up models, and expanding data volume merely optimize the fitting accuracy within the paradigm; they cannot circumvent the structural limitations inherent to the OD paradigm, making it difficult to achieve genuine scientific intelligence and general generalization. Consequently, it is of significant theoretical and applied value to re-examine the inherent trade-offs of the OD paradigm and to explore entirely new solving paradigms tailored for scientific problems.

\section{Ontological Limitations of the Optimization-Driven (OD) Paradigm}
The Optimization-Driven paradigm simplifies complex real-world scientific problems into search problems within a sample data space, seeking an overall optimal solution through globally weighted balancing of all sample constraints. However, real observational data universally contain noise, redundancy, and contradictions. The corresponding solving system is essentially an overdetermined system of conflicting constraints, which gives rise to two inherent limitations of the paradigm:

\subsection{Statistical Compromise under Conflicting Constraints}
The optimization logic of global searching naturally necessitates the model to compromise among various conflicting constraints, ultimately generating a mediocre, compromise solution that fits the overall statistical distribution. In the gradient descent process devoid of structural constraints, the algorithm produces a strong smoothing effect, actively weakening and flattening sparse, localized, and fine-grained features, as well as singular structures within the training set.

This mechanism directly causes the model to overfit statistical patterns of data rather than objective physical laws. Guided by the task evaluation system, the model is prone to ``Reward Hacking'' behaviors targeting the evaluation function \cite{amodei2016}, generating false results with excellent metrics but violating true mechanisms. This creates a scientific illusion of ``optimal metrics, distorted mechanisms,'' failing to accurately restore niche, fine-grained, and critical objective structural features \cite{hochberg2026}.

\subsection{Blind Iterative Optimization without Mechanistic Constraints}
The parameter iteration process of traditional deep learning, whose core engine is the BP algorithm, traces back to the foundational work of Rumelhart, Hinton, and Williams \cite{rumelhart1986}. However, since its inception, this algorithm has lacked constraints from physical laws, mathematical structures, and causal mechanisms; its essence remains a trial-and-error blind search within a high-dimensional non-convex parameter space. Whether it involves optimizing iteration strategies, refining the BP algorithm, or expanding computing power and data volume, all are merely adjustments within the paradigm and cannot alter the underlying logic of ``statistical fitting and blind searching.''

This paradigm can only learn correlated features within the training data distribution and lacks the ability for autonomous deduction and mechanistic extrapolation. Once facing unknown scenarios, novel structures, and boundary conditions in Out-of-Distribution (OOD) environments, the model's fitting system loses its statistical foundation, and reasoning performance decays significantly, resulting in a fundamental generalization bottleneck \cite{tikhonov1977,engl1996}.

\section{Core Reconstruction: Structure-Enabled Intelligence to Break the Statistical Fitting Bottleneck}
Artificial intelligence is a high-performance computing tool; its intelligence originates from rules and logic endowed by humans, not from autonomous model cognition. The limitation of the OD paradigm is not a defect in computing power, data, or algorithm iteration, but a misalignment in the dominant logic of problem-solving. The core pathway to breaking the existing bottleneck is to step out of the inherent framework of data-driven objectives and establish a new Structure-Driven (SD), mechanism-dominated solving paradigm, thereby realizing ``Structure-Enabled Intelligence'' for AI.

\subsection{Constructing Foundational Algorithmic Constraints with Objective Mechanisms}
By abandoning the logic of ``data-dominant solving and global error priority,'' objective physical and mathematical structures---such as conservation laws, wave equations, boundary conditions, and mathematical operators---are embedded into the underlying algorithmic architecture as rigid constraints for model solving. This ensures that the model's iterative calculation, feature extraction, and result generation fully comply with the inevitability of objective mechanisms, constraining the solving domain within a compact set that conforms to physical laws. This fundamentally suppresses the statistical smoothing effect and eliminates fitting defects that distort real structures to adapt to data distributions \cite{tikhonov1977,engl1996}.

This shift is not an isolated engineering improvement. Behind it lies an emerging cross-domain methodological logic: when facing complex problems, the primary question is no longer ``how to solve it faster,'' but ``what is the structure of the problem'' \cite{chen2026sd}.

\subsection{Replacing Black-Box Statistical Training with Structured Solutions}
Relying on structured methods such as physical adjoint theory, virtual source projection, and operator constraints, this approach replaces the traditional pure black-box chain-rule backpropagation mechanism. For example, the Structure-Driven Inversion (SDI) paradigm utilizes algebraic structures such as unitary diagonalization of matrices to achieve the Direct Construction of Operators. Its Physical-Structure-Driven Backpropagation (PSDBP) method achieves structural projection through analogy with physical systems, providing a transparent solution path for deep neural networks \cite{chen2026sdi}. Furthermore, operator learning networks based on functional approximation (such as physics-informed DeepONets) \cite{lu2021} jointly construct a mechanism-transparent and structurally controllable direct solving system, weakening the dependence on massive iterative data training.

Under the Structure-Enabled Intelligence paradigm, the role of data undergoes a fundamental transformation: it no longer serves as the core basis for training fitting and trade-off compromises, but is utilized solely to verify, correct, and supplement structured solving results. The model's core capability originates from built-in scientific structures rather than data statistics, thereby obtaining stable OOD generalization capabilities and interpretability.

\section{Reconstructing Human-AI Synergy: Upholding Human Subjectivity in Scientific Cognition}
Computing-power-driven AI cannot autonomously generate scientific cognition, propose scientific questions, or establish causal systems. It also lacks human rational judgment, innovative thinking, and a sense of responsibility. The current mainstream ``AI for Science'' research mostly remains at the superficial level of data fitting empowering scientific computing, without achieving true mechanistic empowerment and intelligent innovation.

A higher-order scientific intelligence pathway is ``Science for AI'': utilizing mature scientific theories, mathematical structures, and physical mechanisms to discipline artificial intelligence, and using human scientific cognition to define the model's solving logic and boundary rules. This pathway is by no means empty talk. The Structure-Driven Methodology (SD) \cite{chen2026sd} and the Structure-Driven Inversion (SDI) paradigm \cite{chen2026sdi} (particularly its PSDBP algorithmic implementation) are concrete practices of this pathway at the mathematical and engineering levels. In this new paradigm system, the irreplaceable responsibility of mathematicians and scientists lies in the structural definition and Well-posedness analysis of functional spaces \cite{tikhonov1977}.

In the new paradigm system, the core of human-AI synergy is no longer the machine's computing power advantage, but the human cognitive advantage. The core value of human scholars does not lie in computing and fitting, but in defining structures, implanting mechanisms, proposing questions, discerning authenticity, and bearing decision-making responsibilities. Only by upholding the human subject position of intelligence-enabling, and utilizing scientific structures to master computational tools, can we avoid blind reliance on algorithmic statistical results and safeguard the rational essence of scientific research.

\section{Conclusion}
Computing power scaling, data fitting, and iterative optimization can only achieve optimal fitting at the statistical level, incapable of replacing human scientific insight and mechanistic cognition. The statistical compromise characteristic of the OD paradigm is the underlying root cause of insufficient AI generalization, mechanistic distortion, and poor interpretability; all improvements within the paradigm encounter an insurmountable ceiling.

From objective compromise to structure-enabled intelligence, and from data-driven to mechanism-dominant, is an inevitable paradigm leap in the development of scientific intelligence. Abandoning the inherent path of black-box statistical fitting and reconstructing the algorithmic foundation with objective mathematical and physical structures---endowing models with mechanistic logic and extrapolation capabilities through human intelligence-enabling---can break the current generalization bottleneck of AI. It can also reshape the human-AI collaborative scientific research model, enabling artificial intelligence to truly serve the exploration and restoration of objective laws.

\end{document}